\documentclass[12pt]{article}
\usepackage{amssymb,amsmath,amsthm,amsfonts,enumerate, bbm, mathdots}

\usepackage{latexsym}
\usepackage{amscd}
\usepackage{stmaryrd}
\usepackage{color}
\usepackage[all]{xypic}
\usepackage{epsfig}
\usepackage{graphics}
\usepackage{ifthen}
\usepackage{varioref}
\usepackage{rotating}
\usepackage{extarrows}
\usepackage{cite}
\usepackage{mathrsfs}

\numberwithin{equation}{section}

\theoremstyle{plain}
\newtheorem{thm}{Theorem}[section]
\newtheorem{cor}[thm]{Corollary}
\newtheorem{lem}[thm]{Lemma}

\definecolor{darkgreen}{rgb}{0.0625,0.64,0.0625}
\usepackage[%dvipdfm,
pdfstartview=FitH,
CJKbookmarks=true,
bookmarksnumbered=true,
bookmarksopen=true,
colorlinks,
pdfborder=001,
linkcolor=blue,
anchorcolor=green,
citecolor=red
]{hyperref}

\newfont{\scyr}{wncyr10 scaled 550}

\def\proof{\noindent {\bf Proof.\;}}

\allowdisplaybreaks

\begin{document}
	
\title{Generating functions of multiple $t$-star values\thanks{This work was supported by the Fundamental Research Funds for the Central Universities (grant number 22120210552).}}
	
\date{\small ~ \qquad\qquad School of Mathematical Sciences, Tongji University \newline No. 1239 Siping Road, Shanghai 200092, China}
	
\author{Zhonghua Li\thanks{Email address: zhonghua\_li@tongji.edu.cn} ~and ~Lu Yan\thanks{Email address: 1910737@tongji.edu.cn}}

\maketitle
	
\begin{abstract}
In this paper, we study the generating functions of multiple $t$-star values with an arbitrary number of blocks of twos, which are based on the results of the corresponding generating functions of multiple $t$-harmonic star sums. These generating functions can deduce an explicit expression of multiple $t$-star values. As applications, we obtain some evaluations of multiple $t$-star values with one-two-three or more general indices.
\end{abstract}
	
	
	
%%---------------------------------------------------------------------------
%%------------------------Content-------------------------------------------
%%----------------------------------------------------------------------------
	
\section{Introduction}\label{Sec:Intro}
	
Let $\mathbb{C}$, $\mathbb{R}$ and $\mathbb{N}$ be the set of complex numbers, real numbers and positive integers, respectively. And let $\mathbb{N}_0=\mathbb{N}\cup \{0\}$. A finite sequence  $\boldsymbol{s}=(s_1,\ldots,s_r)$ of positive integers is called an index. If $s_1>1$, the multiple zeta value $\zeta(\boldsymbol{s})$ and the multiple zeta-star value $\zeta^{\star}(\boldsymbol{s})$ are defined respectively by (cf. \cite{Hoffman1992,DZagier1994})
\begin{align*}
\zeta(\boldsymbol{s})\equiv \zeta(s_1,\ldots,s_r):=\sum_{k_1>\cdots>k_r\geq 1}
\frac{1}{k_1^{s_1}\cdots k_r^{s_r}},\\
\zeta^{\star}(\boldsymbol{s})\equiv \zeta^{\star}(s_1,\ldots,s_r):=\sum_{k_1\geq\cdots\geq k_r\geq 1}
\frac{1}{k_1^{s_1}\cdots k_r^{s_r}}.
\end{align*}
In 2012, D. Zagier \cite{DZagier2012} obtained explicit evaluation formulas for $\zeta(\{2\}^a,3,\{2\}^b)$ and $\zeta^{\star}(\{2\}^{a},3,\{2\}^{b})$ by establishing the generating functions, where $a,b\in\mathbb{N}_0$ and the $\{s\}^a$ denotes the sequence obtained by repeating $s$ exactly $a$ times. After that, other proofs of Zagier's evaluation formulas have appeared. For example, see \cite{Lee-Peng2018,Li2013,Li-Lupu-Orr} for the proof of the formula for $\zeta(\{2\}^a,3,\{2\}^b)$, and \cite{PPT2014} for the proof of the formula for $\zeta^{\star}(\{2\}^{a},3,\{2\}^{b})$. Especially in 2017, Kh. Hessami Pilehrood and T. Hessami Pilehrood \cite{PP2017} gave another proof of the theorem of Zagier by using a new representation of the generating function of $\zeta^{\star}(\{2\}^{a},3,\{2\}^{b})$. Later in \cite{PP2019}, Kh. Hessami Pilehrood and T. Hessami Pilehrood studied the generating function of multiple zeta-star values with an arbitrary number of blocks of twos. As applications, they \cite{PP2022} provided another way to derive some known evaluations of multiple zeta-star values, and gave new explicit evaluations of $\zeta^{\star}(\{\{2\}^a,3,\{2\}^a,1\}^d)$ and $\zeta^{\star}(\{\{2\}^a,3,\{2\}^a,1\}^d,\{2\}^{a+1})$, where $a\in\mathbb{N}_0$ and $d\in\mathbb{N}$. Note that there have been numerous contributions on evaluations of multiple zeta values and multiple zeta-star values, for example, see \cite{BBBL1998,BBBL2001, PPT2014, SMuneta2008, Ohno-Zudilin, DZagier1994, Zhao2016} and the references therein.
	
In this paper, we consider the evaluations of multiple $t$-star values. Recall from \cite{Hoffman2019} that for an index $\boldsymbol{s}=(s_1,\ldots,s_r)$ with $s_1>1$, Hoffman's multiple $t$-value $t(\boldsymbol{s})$ and multiple $t$-star value $t^{\star}(\boldsymbol{s})$ are defined by
\begin{align*}
t(\boldsymbol{s})&\equiv t(s_1,\ldots,s_r):=\sum_{k_1>\cdots>k_r\geq 1\atop k_i:\text{odd}} \frac{1}{k_1^{s_1}\cdots k_r^{s_r}}\\
&=\sum_{k_1>\cdots>k_r\geq 1} \frac{1}{(2k_1-1)^{s_1}\cdots (2k_r-1)^{s_r}},\\
t^{\star}(\boldsymbol{s})&\equiv t^{\star}(s_1,\ldots,s_r):=\sum_{k_1\geq\cdots\geq k_r\geq 1\atop k_i: \text{odd}} \frac{1}{k_1^{s_1}\cdots k_r^{s_r}}\\
&=\sum_{k_1\geq\cdots\geq k_r\geq 1} \frac{1}{(2k_1-1)^{s_1}\cdots (2k_r-1)^{s_r}},
\end{align*}
respectively. Here if the index is empty, we treat the values $t(\emptyset)=t^{\star}(\emptyset)=1$. By the method of establishing the generating functions, T. Murakami\cite{TMurakami} obtained the evaluation of $t(\{2\}^a,3,\{2\}^b)$ and S. Charlton\cite{SCharlton} gave the evaluation of $t(\{2\}^a,1,\{2\}^b)$. In a recent paper\cite{Li-Lupu-Orr}, L. Lai, C. Lupu and D. Orr provided another proof of Murakami's result of $t(\{2\}^a,3,\{2\}^b)$. In \cite{Li-Wang}, Z. Li and Z. Wang found evaluations of $t^{\star}(\{2\}^a,3,\{2\}^b)$ and $t^{\star}(\{2\}^a,1,\{2\}^b)$.
	
In this paper, we study the generating functions of multiple $t$-star values with an arbitrary number of blocks of twos. Inspired by \cite{PP2019}, we first consider the generating functions of multiple $t$-harmonic star sums which is a finite version of multiple $t$-star values. Then we prove that the limit transfer from the generating functions of multiple $t$-harmonic star sums to the generating functions of multiple $t$-star values is correct. Therefore, we can obtain explicit expressions of the generating functions of multiple $t$-star values.  As applications, we get similar evaluation formulas as in \cite{Li-Wang} for $t^{\star}(\{2\}^a,3,\{2\}^b)$ and $t^{\star}(\{2\}^a,1,\{2\}^b)$ by another way, and show that the multiple $t$-star values with more general index can be evaluated by the weighted sum formulas of alternating multiple $t$-values. Recall from \cite{Quan2020} that for $\boldsymbol{s}=(s_1,\ldots,s_r)\in\mathbb{N}^r$ and $\boldsymbol{\sigma}=(\sigma_1,\ldots,\sigma_r)\in\{\pm1\}^r$ with $(s_1,\sigma_1)\neq (1,1)$, the alternating multiple $t$-values are defined by
\begin{align*}
t(\boldsymbol{s};\boldsymbol{\sigma})\equiv t(s_1,\ldots,s_r;\sigma_1,\ldots,\sigma_r)
:=\sum_{k_1>\cdots>k_r\geq 1}\frac{\sigma_1^{k_1}\cdots\sigma_r^{k_r}}
{(2k_1-1)^{s_1}\cdots(2k_r-1)^{s_r}}.
\end{align*}
	
This paper is organized as follows. In Section \ref{Sec: Generating functions of multiple $t$-harmonic star sums}, we introduce the multiple $t$-harmonic star sums and give their generating functions. In Section \ref{Sec: Multiple $t$-star values and their generating functions}, we obtain the generating functions of multiple $t$-star values and deduce an explicit expression of multiple $t$-star values with an arbitrary number of blocks of twos. As applications, we get some evaluations of multiple $t$-star values with one-two-three or more general indices in Section \ref{Sec: Evaluations for MtSVs}.
		
The following notations will be used throughout the paper. For any index $\boldsymbol{s}=(s_1,\ldots,s_r)$ and $k, m \in\mathbb{N}_0$, we define that
\begin{equation}\label{definitionofV}
V_{k,m}^{\#}(\boldsymbol{s})=
\begin{cases}
\frac{2k-1}{2m-1}\sum_{k\geq l_1\geq\cdots\geq l_r\geq m}\frac{2^{\triangle(k,l_1)+\triangle(l_1,l_2)+\cdots+\triangle(l_r,m)}}{(2l_1-1)^{s_1}\cdots(2l_r-1)^{s_r}} &\text{if\;} \boldsymbol{s}\neq\emptyset \text{\;and\;}k\geq m,\\
&\\
\left(2\cdot\frac{2k-1}{2m-1}\right)^{\triangle(k,m)} &\text{otherwise},
\end{cases}	
\end{equation}
where
\[
\triangle(k,m)=
\begin{cases}
		0 &\text{if\;} k=m,\\
		1 &\text{if\;} k\neq m.
\end{cases}
\]
Also, for any $c\in\mathbb{N}_0$, define
$$\delta(c)=
\begin{cases}
2 &\text{if\;} c=0,\\
1 &\text{if\;} c=1,\\
0 &\text{if\;} c\geq3.
\end{cases}$$

%%%%%%%%%%%%%%%%%%%%%%%%%%%%%%%%%%%%%%%%%%%%%%%%%%%%%%%%%%%%%%%%%%%%%%%%%%%%%%%%%%%%%%%%%%%%%%%%%%%%%%%%%%%%%%%%%
	
\section{Generating functions of multiple $t$-harmonic star sums}\label{Sec: Generating functions of multiple $t$-harmonic star sums}
	
In this section, let $n$ be a fixed positive integer. For any index $\boldsymbol{s}=(s_1,\ldots,s_r)$, we define the multiple $t$-harmonic star sum $t^{\star}_n(\boldsymbol{s})$ by
\begin{align*}
t^{\star}_n(\boldsymbol{s}):= \sum_{n\geq k_1\geq\cdots\geq k_r\geq1} \frac{1}{ (2k_1-1)^{s_1}\cdots(2k_r-1)^{s_r}}.
\end{align*}
Here we set $t^{\star}_n(\emptyset)=1$. For $d\in\mathbb{N}_0$, $\boldsymbol{c}=(c_1,\ldots,c_d)\in\mathbb{N}^d$ and $\boldsymbol{z}=(z_0,z_1,\ldots,z_d)\in\mathbb{C}^{d+1}$, we denote the generating function of multiple $t$-harmonic star sums with an arbitrary number of blocks of twos as
\begin{align*}
		G_n(\boldsymbol{c};\boldsymbol{z})=\sum_{a_0,a_1,\ldots,a_d\geq0}	t^{\star}_n(\{2\}^{a_0},c_1,\{2\}^{a_1},\ldots,c_d,\{2\}^{a_d})z_0^{2a_0}z_1^{2a_1}\cdots z_d^{2a_d},
\end{align*}
Recall that the $\{s\}^a$ denotes the sequence obtained by repeating $s$ exactly $a$ times.

The following  theorem gives an explicit expression of the generating function of multiple $t$-harmonic star sums.
	
\begin{thm}\label{G_n-main-thm}
For any $n\in\mathbb{N}$, $d\in\mathbb{N}_0$, $\boldsymbol{c}=(c_1,\ldots,c_d)\in(\mathbb{N}\setminus\{2\})^d$ and $\boldsymbol{z}=(z_0,z_1,\ldots,z_d)\in\mathbb{C}^{d+1}$ with $|z_j|<1, j=0,1,\ldots,d$, we have
\begin{align}\label{G_n-main-formula}
G_n(\boldsymbol{c};\boldsymbol{z})=\frac{n\binom{2n}{n}}{2^{4n-2}}\sum_{n\geq k_0\geq k_1\geq\cdots\geq k_d\geq1}\binom{2n-1}{n-k_0}\prod_{i=0}^{d}\frac{(-1)^{k_i\delta_i}(2k_i-1)^{\delta_i-1}}{(2k_i-1)^2-z_i^2}V_{k_{i-1},k_i}^{\#}(\{1\}^{c_i-3}),
\end{align}
where $k_{-1}=0$, $\delta_i=\delta(c_i)+\delta(c_{i+1})$ with $c_0=1$ and $c_{d+1}=0$.
\end{thm}

We first give some corollaries of Theorem \ref{G_n-main-thm} in Subsection \ref{SubSec:CorThm2-1}. Then we prepare some preliminary lemmas in Subsection \ref{SubSec:PreLemThm2-1}. Finally, we prove Theorem \ref{G_n-main-thm} in Subsection \ref{SubSec:ProfThm2-1}.
	
\subsection{Corollaries of Theorem \ref{G_n-main-thm}} \label{SubSec:CorThm2-1}
	
Taking $\boldsymbol{c}=(\{3,1\}^d)$ in \eqref{G_n-main-formula}, we get the following corollary.
	
\begin{cor}
For any $n\in\mathbb{N}$, $d\in\mathbb{N}_0$ and any complex numbers $z_0,z_1,\ldots,z_{2d}$ with $|z_j|<1, j=0,1,\ldots,2d$, we have
\begin{align*}
&\sum_{a_0,a_1,\ldots,a_{2d}\geq0}t^{\star}_n(\{2\}^{a_0},3,\{2\}^{a_1},1,\{2\}^{a_2},\ldots,3,\{2\}^{a_{2d-1}},1,\{2\}^{a_{2d}})z_0^{2a_0}z_1^{2a_1}\cdots z_{2d}^{2a_{2d}}\\
=&\frac{n\binom{2n}{n}}{2^{4n-2}}\sum_{n\geq k_0\geq\cdots\geq k_{2d}\geq1}\binom{2n-1}{n-k_0}(2k_{2d}-1)^2\prod_{i=0}^{2d}\frac{(-1)^{k_i}}{(2k_i-1)^2-z_i^2}\left(2\cdot\frac{2k_{i-1}-1}{2k_{i}-1}\right)^{\triangle(k_{i-1},k_i)},
\end{align*}
where $k_{-1}=0$.
\end{cor}
	
Removing the case of $z_u=0$ for $0\leq u\leq d$, we have
	
\begin{cor}\label{a_u>=1-corollary}
Let $n\in\mathbb{N}$, $d\in\mathbb{N}_0$, $c_1,\ldots,c_d\in\mathbb{N}\setminus\{2\}$ and $z_0,z_1,\ldots,z_d$ be complex numbers with $|z_j|<1, j=0,1,\ldots,d$. Then for any integer $u$ with $0\leq u\leq d$, we have
\begin{align*}
&\sum_{a_0,a_1,\ldots,a_d\geq0\atop a_u\geq1}t^{\star}_n(\{2\}^{a_0},c_1,\{2\}^{a_1},\ldots,c_d,\{2\}^{a_d})z_0^{2a_0}z_1^{2a_1}\cdots z_d^{2a_d}\\
=&z_u^2\cdot\frac{n\binom{2n}{n}}{2^{4n-2}}\sum_{n\geq k_0\geq\cdots\geq k_d\geq1}\binom{2n-1}{n-k_0}\frac{1}{(2k_u-1)^2}\prod_{i=0}^{d}\frac{(-1)^{k_i\delta_i}(2k_i-1)^{\delta_i-1}}{(2k_i-1)^2-z_i^2}V_{k_{i-1},k_i}^{\#}(\{1\}^{c_i-3}),
\end{align*}
where $k_{-1}=0$, $\delta_i=\delta(c_i)+\delta(c_{i+1})$ with $c_0=1$ and $c_{d+1}=0$.
\end{cor}
	
\proof
We observe that
\begin{align}\label{G_n-G_n with z_u=0}
&\sum_{a_0,a_1,\ldots,a_d\geq0\atop a_u\geq1}t^{\star}_n(\{2\}^{a_0},c_1,\{2\}^{a_1},\ldots,c_d,\{2\}^{a_d})z_0^{2a_0}z_1^{2a_1}\cdots z_d^{2a_d}\nonumber\\
=&G_n(\boldsymbol{c}; \boldsymbol{z})-G_n(\boldsymbol{c}; z_0,\ldots,z_{u-1}, 0, z_{u+1},\ldots,z_d),
\end{align}
where $\boldsymbol{c}=(c_1,\ldots,c_d)$ and $\boldsymbol{z}=(z_0,z_1,\ldots,z_d)$. Hence we get the desired result by Theorem \ref{G_n-main-thm}.
\qed
	
Expanding the right-hand side of \eqref{G_n-main-formula}, we obtain the following result.
	
\begin{cor}\label{t_n^star-corollary}
For any $n\in\mathbb{N}$, $d\in\mathbb{N}_0$, $c_1,\ldots,c_d\in\mathbb{N}\setminus\{2\}$ and    $a_0,a_1,\ldots,a_d\in\mathbb{N}_0$, we have
\begin{align*}
&t^{\star}_n(\{2\}^{a_0},c_1,\{2\}^{a_1},\ldots,c_d,\{2\}^{a_d})\\
=&\frac{n\binom{2n}{n}}{2^{4n-2}}\sum_{n\geq k_0\geq k_1\geq\cdots\geq k_d\geq1}\binom{2n-1}{n-k_0}\prod_{i=0}^{d}\frac{(-1)^{k_i\delta_i}}{(2k_i-1)^{2a_i-\delta_i+3}}V_{k_{i-1},k_i}^{\#}(\{1\}^{c_i-3}),
\end{align*}
where $k_{-1}=0$, $\delta_i=\delta(c_i)+\delta(c_{i+1})$ with $c_0=1$ and $c_{d+1}=0$.
\end{cor}
	
\proof
Applying the power series expansion
$$\frac{1}{(2k_i-1)^2-z_i^2}=\sum\limits_{a_i=0}^\infty\frac{z_i^{2a_i}}{(2k_i-1)^{2a_i+2}}$$
for $i=0,1,\ldots,d$ and comparing the coefficients of $z_0^{2a_0}z_1^{2a_1}\cdots z_d^{2a_d}$ in \eqref{G_n-main-formula}, we get the result.
\qed

\subsection{Preliminary Lemmas}\label{SubSec:PreLemThm2-1}
	
Similar to \cite[(2.5)]{PP2019} and \cite[(2.1), (2.2)]{PPT2014}, we prove several finite sum formulas that will be used in the proof of Theorem \ref{G_n-main-thm}.
	
\begin{lem}
For any $n\in\mathbb{N}$ and $l\in\mathbb{N}_0$, we have
\begin{align}\label{l+1-n}
\sum_{k=l+1}^n(2k-1)\binom{2n-1}{n-k}=(n-l)\binom{2n-1}{n-l}
\end{align}
and
\begin{align}\label{alter-l+1-n}
\sum_{k=l+1}^n(-1)^k\binom{2n-1}{n-k}=\frac{(-1)^{l+1}(n-l)}{2n-1}\binom{2n-1}{n-l}.
\end{align}
\end{lem}
	
\proof
We may assume that $n\geq l+1$. For any integer $k$ with $l+1\leq k\leq n$, we observe that
\begin{align}\label{I-recursion}
(2k-1)\binom{2n-1}{n-k}=I(n;k+1)-I(n;k)
\end{align}
with
\begin{align*}
I(n;j)=-(n+j-1)\binom{2n-1}{n-j}
\end{align*}
for $l+1\leq j\leq n$	and $I(n;n+1)=0$. Summing both sides of \eqref{I-recursion} over $k$ from $l+1$ to $n$, \eqref{l+1-n} follows easily.
	
Similarly, to prove \eqref{alter-l+1-n}, we set
\begin{align*}
J(n;j)=\frac{(-1)^{j-1}(n+j-1)}{2n-1}\binom{2n-1}{n-j}
\end{align*}
for $l+1\leq j\leq n$	and $J(n;n+1)=0$. Then it is easy to see that for any integer $k$ with $l+1\leq k\leq n$,
\begin{align*}
(-1)^k\binom{2n-1}{n-k}=J(n;k+1)-J(n;k).
\end{align*}
Then summing both sides of the above equation over $k$ from $l+1$ to $n$, we obtain \eqref{alter-l+1-n}.
\qed
	
\begin{lem}\label{alter-l-n-V}
For any $n,l\in\mathbb{N}$ and $c\in\mathbb{N}_0$, we have
\begin{align}\label{l-n-V}
\sum_{k=l}^n\frac{(-1)^k}{2k-1}\binom{2n-1}{n-k}V_{k,l}^{\#}(\{1\}^c)=\frac{(-1)^l}{(2n-1)^{c+1}}\binom{2n-1}{n-l}.
\end{align}
\end{lem}
	
\proof
We may assume that $n\geq l$. The proof is by induction on $c$. For $c=0$, it is enough to show
\begin{align}\label{c=0}
\sum_{k=l}^n\frac{(-1)^k}{2k-1}\binom{2n-1}{n-k}\left(2\cdot\frac{2k-1}{2l-1}\right)^{\triangle(k,l)}
=\frac{(-1)^l}{2n-1}\binom{2n-1}{n-l}.
\end{align}
In fact, the left-hand side of \eqref{c=0} is
$$\frac{(-1)^l}{2l-1}\binom{2n-1}{n-l}+\frac{2}{2l-1}\sum_{k=l+1}^n(-1)^k\binom{2n-1}{n-k}.$$
Applying \eqref{alter-l+1-n}, we find the left-hand side of \eqref{c=0} becomes
$$\frac{(-1)^l}{2l-1}\binom{2n-1}{n-l}+\frac{2(-1)^{l+1}(n-l)}{(2l-1)(2n-1)}\binom{2n-1}{n-l},$$
which exactly is $\frac{(-1)^l}{2n-1}\binom{2n-1}{n-l}$.
	
If $c\geq1$, the left-hand side of \eqref{l-n-V} is
$$\sum_{k=l}^n\frac{(-1)^k}{2k-1}\binom{2n-1}{n-k}\cdot\frac{2k-1}{2l-1}\sum_{k\geq l_1\geq\cdots\geq l_c\geq l}\frac{2^{\triangle(k,l_1)+\triangle(l_1,l_2)+\cdots+\triangle(l_c,l)}}{(2l_1-1)\cdots(2l_c-1)}.$$
Changing the order of the summations, the left-hand side of \eqref{l-n-V} becomes
$$\sum_{n\geq l_1\geq\cdots\geq l_c\geq l}\left\{\sum_{k=l_1}^{n}\frac{(-1)^k2^{\triangle(k,l_1)}}{2l_1-1}\binom{2n-1}{n-k}\right\}
\frac{2^{\triangle(l_1,l_2)+\cdots+\triangle(l_c,l)}}{(2l-1)(2l_2-1)\cdots(2l_c-1)},$$
which is
$$\sum_{n\geq l_1\geq\cdots\geq l_c\geq l}\left\{\sum_{k=l_1}^{n}\frac{(-1)^k}{2k-1}\binom{2n-1}{n-k}\left(2\cdot\frac{2k-1}{2l_1-1}\right)^{\triangle(k,l_1)}\right\}\frac{2^{\triangle(l_1,l_2)+\cdots+\triangle(l_c,l)}}{(2l-1)(2l_2-1)\cdots(2l_c-1)}.$$
Using \eqref{c=0} to deal with the part in the curly brace, we find the left-hand side of \eqref{l-n-V} is
$$\sum_{n\geq l_1\geq\cdots\geq l_c\geq l}\frac{(-1)^{l_1}}{2n-1}\binom{2n-1}{n-l_1}\cdot\frac{2^{\triangle(l_1,l_2)+\cdots+\triangle(l_c,l)}}{(2l-1)(2l_2-1)\cdots(2l_c-1)},$$
which is just
$$\frac{1}{2n-1}\sum_{l_1=l}^n\frac{(-1)^{l_1}}{2l_1-1}\binom{2n-1}{n-l_1}V_{l_1,l}^{\#}(\{1\}^{c-1}).$$
Then the result follows from the inductive hypothesis for $c-1$.
\qed

\subsection{Proof of Theorem \ref{G_n-main-thm}} \label{SubSec:ProfThm2-1}

We prove Theorem \ref{G_n-main-thm} by induction on $n+d$. When $n=1$, we have
\begin{align*}
G_1(\boldsymbol{c}; \boldsymbol{z})=\sum_{a_0,a_1,\ldots,a_d\geq0}z_0^{2a_0}\cdots z_d^{2a_d}=\prod_{i=0}^{d}\frac{1}{1-z_i^2},
\end{align*}
and the right-hand side of \eqref{G_n-main-formula} is
\begin{align*}
\frac{1}{2}\prod_{i=0}^d\frac{(-1)^{\delta_i}}{1-z_i^2}\cdot (-2)=(-1)^{1+\delta_0+\delta_1+\cdots+\delta_d}\prod_{i=0}^d\frac{1}{1-z_i^2}=\prod_{i=0}^d\frac{1}{1-z_i^2}.
\end{align*}
Hence \eqref{G_n-main-formula} is true for $n=1$.
	
When $d=0$, since for $n>1$,
$$t_n^{\star}(\{2\}^{a_0})=\sum_{k=0}^{a_0}\frac{t_{n-1}^{\star}(\{2\}^k)}{(2n-1)^{2(a_0-k)}},$$
we have a recursive formula
$$G_n( ; z_0)=\frac{(2n-1)^2}{(2n-1)^2-z_0^2}G_{n-1}( ; z_0).$$
Using the above formula repeatedly, we get
\begin{align*}
G_n( ; z_0)=\prod_{j=2}^{n}\frac{(2j-1)^2}{(2j-1)^2-z_0^2}G_{1}( ; z_0)
=[(2n-1)!!]^{2}\prod_{j=1}^{n}\frac{1}{(2j-1)^2-z_0^2}.
\end{align*}
By the method of partial fractional decomposition, we have
\begin{align*}
G_n( ; z_0)
=&\sum_{k=1}^{n}\frac{2(-1)^{k-1}(2k-1)}{(2n-2k)!!(2n+2k-2)!!}\cdot\frac{[(2n-1)!!]^2}{(2k-1)^2-z_0^2}\\
=&\frac{n\binom{2n}{n}}{2^{4n-2}}\sum_{k=1}^{n}\binom{2n-1}{n-k}\frac{(-1)^k(2k-1)^2}{(2k-1)^2-z_0^2}\cdot\frac{-2}{2k-1}.
\end{align*}
So \eqref{G_n-main-formula} is proved for $d=0$.
	
Now assume that $n>1$ and $d>0$. Since
\begin{align*}
t^{\star}_n(\{2\}^{a_0},c_1,\{2\}^{a_1},\ldots,c_d,\{2\}^{a_d})
=&\sum_{k=0}^{a_0}\frac{1}{(2n-1)^{2a_0-2k}}t^{\star}_{n-1}(\{2\}^{k},c_1,\{2\}^{a_1},\ldots,c_d,\{2\}^{a_d})\\
&+\frac{1}{(2n-1)^{2a_0+c_1}}t^{\star}_n(\{2\}^{a_1},c_2,\{2\}^{a_2},\ldots,c_d,\{2\}^{a_d}),
\end{align*}
we find that
\begin{align*}
G_n(\boldsymbol{c}; \boldsymbol{z})
=&\sum_{a_0,\ldots,a_d\geq0}\sum_{k=0}^{a_0}\frac{1}{(2n-1)^{2a_0-2k}}	t^{\star}_{n-1}(\{2\}^{k},c_1,\{2\}^{a_1},\ldots,c_d,\{2\}^{a_d})z_0^{2a_0}\cdots z_d^{2a_d},\\
&+\sum_{a_0\geq0}\frac{z_0^{2a_0}}{(2n-1)^{2a_0+c_1}}G_n(\boldsymbol{c}^-; \boldsymbol{z}^-),
\end{align*}
where $\boldsymbol{c}^-=(c_2,\ldots,c_d)$ and $\boldsymbol{z}^-=(z_1,z_2,\ldots,z_d)$. Then it is easy to get the following recursive formula
\begin{align}\label{recurence-G}
G_n(\boldsymbol{c}; \boldsymbol{z})=\frac{(2n-1)^2}{(2n-1)^2-z_0^2}G_{n-1}(\boldsymbol{c}; \boldsymbol{z})+\frac{(2n-1)^{2-c_1}}{(2n-1)^2-z_0^2}G_n(\boldsymbol{c}^-; \boldsymbol{z}^-).
\end{align}
Let $\widetilde{G}_n(\boldsymbol{c};\boldsymbol{z})$ denote the right-hand side of \eqref{G_n-main-formula}. Therefore according to \eqref{recurence-G}, we need to prove
\begin{align}\label{G_n-main-new}
\frac{(2n-1)^2}{(2n-1)^2-z_0^2}G_{n-1}(\boldsymbol{c}; \boldsymbol{z})-\widetilde{G}_n(\boldsymbol{c};\boldsymbol{z})
=-\frac{(2n-1)^{2-c_1}}{(2n-1)^2-z_0^2}G_n(\boldsymbol{c}^-; \boldsymbol{z}^-).
\end{align}
	
By the induction hypothesis for $G_{n-1}(\boldsymbol{c}; \boldsymbol{z})$, the left hand-side of \eqref{G_n-main-new} is
\begin{align*}
&\frac{(2n-1)^2}{(2n-1)^2-z_0^2}\frac{(n-1)\binom{2n-2}{n-1}}{2^{4n-6}}\sum_{n-1\geq k_0\geq\cdots\geq k_d\geq1}\binom{2n-3}{n-k_0-1}\nonumber\\
&\quad\quad\quad\quad\quad\quad\times\prod_{i=0}^{d}\frac{(-1)^{k_i\delta_i}(2k_i-1)^{\delta_i-1}}{(2k_i-1)^2-z_i^2}V_{k_{i-1},k_i}^{\#}(\{1\}^{c_i-3})\\
&-\frac{n\binom{2n}{n}}{2^{4n-2}}\sum_{n\geq k_0\geq\cdots\geq k_d\geq1}\binom{2n-1}{n-k_0}\prod_{i=0}^{d}\frac{(-1)^{k_i\delta_i}(2k_i-1)^{\delta_i-1}}{(2k_i-1)^2-z_i^2}V_{k_{i-1},k_i}^{\#}(\{1\}^{c_i-3}),
\end{align*}
which is
\begin{align*}
&\sum_{n\geq k_0\geq\cdots\geq k_d\geq1}\prod_{i=0}^{d}\frac{(-1)^{k_i\delta_i}(2k_i-1)^{\delta_i-1}}{(2k_i-1)^2-z_i^2}V_{k_{i-1},k_i}^{\#}(\{1\}^{c_i-3})\\
&\quad\quad\quad\times\left\{\frac{(2n-1)^2}{(2n-1)^2-z_0^2}\cdot\frac{(n-1)\binom{2n-2}{n-1}}{2^{4n-6}}\binom{2n-3}{n-k_0-1}
-\frac{n\binom{2n}{n}}{2^{4n-2}}\binom{2n-1}{n-k_0}\right\}.
\end{align*}	
As the part in braces can be simplified as
$$-\frac{n\binom{2n}{n}}{2^{4n-2}}\binom{2n-1}{n-k_0}\cdot\frac{(2k_0-1)^2-z_0^2}{(2n-1)^2-z_0^2},$$
the left-hand side of \eqref{G_n-main-new} becomes
\begin{align}\label{difference}
&\frac{2}{(2n-1)^2-z_0^2}\cdot\frac{n\binom{2n}{n}}{2^{4n-2}}\sum_{n\geq k_0\geq k_1\geq\cdots\geq k_d\geq1}\binom{2n-1}{n-k_0}\frac{(-1)^{k_0(1+\delta(c_1))}(2k_0-1)^{\delta(c_1)}}{2k_0-1}\nonumber\\
&\quad\quad\quad\quad\quad\quad\quad\times\prod_{i=1}^{d}\frac{(-1)^{k_i\delta_i}(2k_i-1)^{\delta_i-1}}{(2k_i-1)^2-z_i^2}V_{k_{i-1},k_i}^{\#}(\{1\}^{c_i-3}).
\end{align}
Let $\sum\nolimits_{k_0}$ denote the inner sum over $k_0$ in \eqref{difference}, that is
\begin{align*}
\sum\nolimits_{k_0}=\sum_{k_0=k_1}^{n}\binom{2n-1}{n-k_0}\frac{(-1)^{k_0(1+\delta(c_1))}(2k_0-1)^{\delta(c_1)}}{2k_0-1}V_{k_0,k_1}^{\#}(\{1\}^{c_1-1}).
\end{align*}
	
If $c_1=1$, $\delta(c_1)=1$, we obtain
\begin{align*}
\sum\nolimits_{k_0}=&\sum_{k_0=k_1}^{n}\binom{2n-1}{n-k_0}\left(2\cdot\frac{2k_0-1}{2k_1-1}\right)^{\triangle(k_0,k_1)}\\
=&\binom{2n-1}{n-k_1}+\frac{2}{2k_1-1}\sum_{k_0=k_1+1}^{n}\binom{2n-1}{n-k_0}(2k_0-1).
\end{align*}
Then by \eqref{l+1-n}, we have
$$\sum\nolimits_{k_0}=\binom{2n-1}{n-k_1}+\frac{2n-2k_1}{2k_1-1}\binom{2n-1}{n-k_1}=\frac{2n-1}{2k_1-1}\binom{2n-1}{n-k_1}.$$
Using the above result together with \eqref{difference}, we find the left-hand side of \eqref{G_n-main-new} equals
\begin{align*}
&-\frac{2n-1}{(2n-1)^2-z_0^2}\cdot\frac{n\binom{2n}{n}}{2^{4n-2}}\sum_{n\geq k_1\geq\cdots\geq k_d\geq1}\binom{2n-1}{n-k_1}\frac{(-1)^{k_1(1+\delta(c_2))}(2k_1-1)^{\delta(c_2)}}{(2k_1-1)^2-z_1^2}\cdot\frac{-2}{2k_1-1}\\
&\quad\quad\quad\quad\quad\quad\quad\times\prod_{i=2}^{d}\frac{(-1)^{k_i\delta_i}(2k_i-1)^{\delta_i-1}}{(2k_i-1)^2-z_i^2}V_{k_{i-1},k_i}^{\#}(\{1\}^{c_i-3}),
\end{align*}
which by the inductive hypothesis for $G_{n}(\boldsymbol{c}^-;\boldsymbol{z}^-)$ is
$$-\frac{2n-1}{(2n-1)^2-z_0^2}G_{n}(\boldsymbol{c}^-;\boldsymbol{z}^-).$$
Hence the theorem is proved in this case.
	
If $c_1\geq3$, $\delta(c_1)=0$, by Lemma \ref{alter-l-n-V}, we get
\begin{align*}
\sum\nolimits_{k_0}=\sum_{k_0=k_1}^{n}\binom{2n-1}{n-k_0}\frac{(-1)^{k_0}}{2k_0-1}V_{k_0,k_1}^{\#}(\{1\}^{c_1-3})=\frac{(-1)^{k_1}}{(2n-1)^{c_1-2}}\binom{2n-1}{n-k_1}.
\end{align*}
Combining this formula with \eqref{difference}, we obtain the left-hand side of \eqref{G_n-main-new} is
\begin{align*}
&-\frac{(2n-1)^{2-c_1}}{(2n-1)^2-z_0^2}\cdot\frac{n\binom{2n}{n}}{2^{4n-2}}\sum_{n\geq k_1\geq\cdots\geq k_d\geq1}\binom{2n-1}{n-k_1}\frac{(-1)^{k_1(1+\delta(c_2))}(2k_1-1)^{\delta(c_2)}}{(2k_1-1)^2-z_1^2}\cdot\frac{-2}{2k_1-1}\\
&\quad\quad\quad\quad\quad\quad\quad\times\prod_{i=2}^{d}\frac{(-1)^{k_i\delta_i}(2k_i-1)^{\delta_i-1}}{(2k_i-1)^2-z_i^2}V_{k_{i-1},k_i}^{\#}(\{1\}^{c_i-3}),
\end{align*}
which by the inductive hypothesis for $G_{n}(\boldsymbol{c}^-;\boldsymbol{z}^-)$ is
$$-\frac{(2n-1)^{2-c_1}}{(2n-1)^2-z_0^2}G_{n}(\boldsymbol{c}^-;\boldsymbol{z}^-).$$
Therefore, we conclude that $G_n(\boldsymbol{c}; \boldsymbol{z})=\widetilde{G}_n(\boldsymbol{c};\boldsymbol{z})$. This completes the proof.
\qed

%%%%%%%%%%%%%%%%%%%%%%%%%%%%%%%%%%%%%%%%%%%%%%%%%%%%%%%%%%%%%%%%%%%%%%%%%%%%%%%%%%%%%%%%%%%%%%%%%%%%%%%%%%%%%%%%%

\section{Multiple $t$-star values and their generating functions}\label{Sec: Multiple $t$-star values and their generating functions}
	
In this section, we consider the generating functions of multiple $t$-star values.
For $d\in\mathbb{N}_0$, $\boldsymbol{c}=(c_1,\ldots,c_d)\in\mathbb{N}^d$ with $c_1>1$ and $\boldsymbol{z}=(z_0,z_1,\ldots,z_d)\in\mathbb{C}^{d+1}$, we define the generating function
\begin{align*}
G(\boldsymbol{c}; \boldsymbol{z})=\sum_{a_0,a_1,\ldots,a_d\geq0}	t^{\star}(\{2\}^{a_0},c_1,\{2\}^{a_1},\ldots,c_d,\{2\}^{a_d})z_0^{2a_0}z_1^{2a_1}\cdots z_d^{2a_d}.
\end{align*}
Let $n$ tend to infinity in Theorem \ref{G_n-main-thm}, we can get the following theorem.
	
\begin{thm}\label{G-main-thm}
For any $d\in\mathbb{N}_0$, $\boldsymbol{c}=(c_1,\ldots,c_d)\in(\mathbb{N}\setminus\{2\})^d$ with $c_1\geq3$, and $\boldsymbol{z}=(z_0,z_1,\ldots,z_d)\in\mathbb{C}^{d+1}$ with $|z_j|<1, j=0,1,\ldots,d$, we have
\begin{align}\label{G-main-formula}
G(\boldsymbol{c};\boldsymbol{z})=\frac{2}{\pi}\sum_{k_0\geq k_1\geq\cdots\geq k_d\geq1}\prod_{i=0}^{d}\frac{(-1)^{k_i\delta_i}(2k_i-1)^{\delta_i-1}}{(2k_i-1)^2-z_i^2}V_{k_{i-1},k_i}^{\#}(\{1\}^{c_i-3}),
\end{align}
where $k_{-1}=0$, $\delta_i=\delta(c_i)+\delta(c_{i+1})$ with $c_0=1$ and $c_{d+1}=0$.
\end{thm}
	
Setting $\boldsymbol{c}=(\{3,1\}^d)$ in \eqref{G-main-formula}, we deduce the following result.
	
\begin{cor}
For any $d\in\mathbb{N}_0$ and $z_0,z_1,\ldots,z_{2d}\in\mathbb{C}$ with $|z_j|<1, j=0,1,\ldots,2d$, we have
\begin{align*}
&\sum_{a_0,a_1,\ldots,a_{2d}\geq0}t^{\star}(\{2\}^{a_0},3,\{2\}^{a_1},1,\{2\}^{a_2},\ldots,3,\{2\}^{a_{2d-1}},1,\{2\}^{a_{2d}})z_0^{2a_0}z_1^{2a_1}\cdots z_{2d}^{2a_{2d}}\\
=&\frac{2}{\pi}\sum_{k_0\geq k_1\geq\cdots\geq k_{2d}\geq1}(2k_{2d}-1)^2\prod_{i=0}^{2d}\frac{(-1)^{k_i}}{(2k_i-1)^2-z_i^2}\left(2\cdot\frac{2k_{i-1}-1}{2k_{i}-1}\right)^{\triangle(k_{i-1},k_i)},
\end{align*}
where $k_{-1}=0$.
\end{cor}

The following theorem gives explicit expressions of multiple $t$-star values with an arbitrary number of blocks of twos.

\begin{thm}\label{t^star-theorem}
For any $d,a_0,a_1,\ldots,a_d\in\mathbb{N}_0$, $c_1,\ldots,c_d\in\mathbb{N}\setminus\{2\}$ with $c_1\geq3$ if $a_0=0$ and $d\geq1$, we have
\begin{align}\label{t^star-formula}
&t^{\star}(\{2\}^{a_0},c_1,\{2\}^{a_1},\ldots,c_d,\{2\}^{a_d})\nonumber\\
=&\frac{2}{\pi}\sum_{k_0\geq k_1\geq\cdots\geq k_d\geq1}\prod_{i=0}^{d}\frac{(-1)^{k_i\delta_i}}{(2k_i-1)^{2a_i-\delta_i+3}}V_{k_{i-1},k_i}^{\#}(\{1\}^{c_i-3}),
\end{align}
where $k_{-1}=0$, $\delta_i=\delta(c_i)+\delta(c_{i+1})$ with $c_0=1$ and $c_{d+1}=0$.
\end{thm}

Corresponding to Corollary \ref{a_u>=1-corollary}, we have the following similar result for multiple $t$-star values.

\begin{thm}\label{a_u>=1-theorem}
Let $d,u\in\mathbb{N}_0$ with $d\geq u$. For any $c_1,\ldots,c_d\in\mathbb{N}\setminus\{2\}$ with $c_1\geq3$ if $u\geq1$, and any $z_0,z_1,\ldots,z_d\in\mathbb{C}$ with $|z_j|<1, j=0,1,\ldots,d$, we have
\begin{align*}
&\sum_{a_0,a_1,\ldots,a_d\geq0\atop a_u\geq1}t^{\star}(\{2\}^{a_0},c_1,\{2\}^{a_1},\ldots,c_d,\{2\}^{a_d})z_0^{2a_0}z_1^{2a_1}\cdots z_d^{2a_d}\\
=&z_u^2\cdot\frac{2}{\pi}\sum_{k_0\geq k_1\geq\cdots\geq k_d\geq1}\frac{1}{(2k_u-1)^2}\prod_{i=0}^{d}\frac{(-1)^{k_i\delta_i}(2k_i-1)^{\delta_i-1}}{(2k_i-1)^2-z_i^2}V_{k_{i-1},k_i}^{\#}(\{1\}^{c_i-3}),
\end{align*}
where $k_{-1}=0$, $\delta_i=\delta(c_i)+\delta(c_{i+1})$ with $c_0=1$ and $c_{d+1}=0$.
\end{thm}

In Subsection \ref{SubSec:LemsThms3-1,3,4}, we prove some lemmas which can determine the feasibility of the limit transfer from the generating functions of multiple $t$-harmonic star sums to the generating functions of multiple $t$-star values. The proofs of Theorems \ref{G-main-thm}, \ref{t^star-theorem} and \ref{a_u>=1-theorem} are given in Subsection \ref{SubSec:ProfThms3-1,3,4}.

\subsection{Lemmas}\label{SubSec:LemsThms3-1,3,4}

Similar to \cite[Lemma 4.1]{PP2019}, we have the following result.

\begin{lem}\label{limit-G_n=G-lemma}
Let $d\in\mathbb{N}_0$, $\boldsymbol{c}=(c_1,\ldots,c_d)\in(\mathbb{N}\setminus\{2\})^d$ with $c_1\geq3$, and let $\boldsymbol{z}=(z_0,z_1,\ldots,z_d)\in\mathbb{C}^{d+1}$
with $|z_j|<1, j=0,1,\ldots,d$. Then
\begin{align}\label{limit-G_n=G}
\lim_{n\rightarrow\infty}G_n(\boldsymbol{c}; \boldsymbol{z})=G(\boldsymbol{c}; \boldsymbol{z}),
\end{align}
and the convergence is uniform in any closed region $E$: $|z_0|\leq u_0<1, |z_1|\leq u_1<1, \ldots, |z_d|\leq u_d<1$.
\end{lem}

\proof
For $|z|<1$, we find that
\begin{align}\label{sum-t^*_n,m}
\prod_{k=m}^{n}\left(1-\frac{z^2}{(2k-1)^2}\right)^{-1}=\sum_{l=0}^{\infty}t_{n,m}^{\star}(\{2\}^l)z^{2l}.
\end{align}
Here we define that, for integers $n\geq m >0$ and $\boldsymbol{s}=(s_1,\ldots, s_r)\in\mathbb{N}^r$,
\begin{align*}
t_{n,m}^{\star}(\boldsymbol{s})=t_{n,m}^{\star}(s_1,\ldots, s_r)=\sum_{n\geq k_1\geq\cdots\geq k_r\geq m} \frac{1}{ (2k_1-1)^{s_1}\cdots(2k_r-1)^{s_r}}.
\end{align*}
We also set $t^{\star}_{n,m}(\emptyset)=1$. If $s_1>1$, we may let $n$ tend to infinity and define
$$t_{\infty,m}^{\star}(\boldsymbol{s})=\sum_{k_1\geq\cdots\geq k_r\geq m} \frac{1}{ (2k_1-1)^{s_1}\cdots(2k_r-1)^{s_r}}.$$

If $d=0$, using \eqref{sum-t^*_n,m}, we have
\begin{align*}
G_n( ;z_0)=\prod_{k=1}^n\left(1-\frac{z_0^2}{(2k-1)^2}\right)^{-1}
\end{align*}
and
\begin{align*}
G( ;z_0)=\prod_{k=1}^\infty\left(1-\frac{z_0^2}{(2k-1)^2}\right)^{-1}.
\end{align*}
Hence the lemma is true for $d=0$.

Now assume that $d\geq1$. Using \eqref{sum-t^*_n,m}, we get
\begin{align*}
&G_n(\boldsymbol{c}; \boldsymbol{z})\\
=&\sum_{n\geq k_1\geq\cdots\geq k_d\geq1}\frac{\prod\limits_{k=k_1}^n\left(1-\frac{z_0^2}{(2k-1)^2}\right)^{-1}\prod\limits_{k=k_2}^{k_1}\left(1-\frac{z_1^2}{(2k-1)^2}\right)^{-1}\cdots\prod\limits_{k=1}^{k_d}\left(1-\frac{z_d^2}{(2k-1)^2}\right)^{-1}}{(2k_1-1)^{c_1}\cdots(2k_d-1)^{c_d}}.
\end{align*}
Setting
\begin{align*}
&G_{n}^{\infty}(\boldsymbol{c}; \boldsymbol{z})\\
=&\sum_{n\geq k_1\geq\cdots\geq k_d\geq1}\frac{\prod\limits_{k=k_1}^{\infty}\left(1-\frac{z_0^2}{(2k-1)^2}\right)^{-1}\prod\limits_{k=k_2}^{k_1}\left(1-\frac{z_1^2}{(2k-1)^2}\right)^{-1}\prod\limits_{k=1}^{k_d}\left(1-\frac{z_d^2}{(2k-1)^2}\right)^{-1}}{(2k_1-1)^{c_1}\cdots(2k_d-1)^{c_d}},
\end{align*}
then
\begin{align}\label{G_n^inf-G_n}
&|G_{n}^{\infty}(\boldsymbol{c}; \boldsymbol{z})-G_{n}(\boldsymbol{c}; \boldsymbol{z})|\nonumber\\
\leq&\sum_{n\geq k_1\geq\cdots\geq k_d\geq1}\frac{\prod\limits_{k=k_1}^{n}\left(1-\frac{|z_0|^2}{(2k-1)^2}\right)^{-1}\prod\limits_{k=k_2}^{k_1}\left(1-\frac{|z_1|^2}{(2k-1)^2}\right)^{-1}\prod\limits_{k=1}^{k_d}\left(1-\frac{|z_d|^2}{(2k-1)^2}\right)^{-1}}{(2k_1-1)^{c_1}\cdots(2k_d-1)^{c_d}}\nonumber\\
&\quad\quad\quad\quad\quad\quad\quad\quad\quad\quad\quad\quad\quad\times\left|\prod\limits_{k=n+1}^{\infty}\left(1-\frac{z_0^2}{(2k-1)^2}\right)^{-1}-1\right|.
\end{align}
Assume that $|z_0|\leq u_0<1, |z_1|\leq u_1<1, \ldots, |z_d|\leq u_d<1$. Using the infinite product formula for the sine function
$$\frac{\sin\pi z}{\pi z}=\prod\limits_{k=1}^\infty\left(1-\frac{z^2}{k^2}\right),$$
we obtain
\begin{align}\label{G_n-upper-bound}
&\sum_{n\geq k_1\geq\cdots\geq k_d\geq1}\frac{\prod\limits_{k=k_1}^n\left(1-\frac{|z_0|^2}{(2k-1)^2}\right)^{-1}\prod\limits_{k=k_2}^{k_1}\left(1-\frac{|z_1|^2}{(2k-1)^2}\right)^{-1}\prod\limits_{k=1}^{k_d}\left(1-\frac{|z_d|^2}{(2k-1)^2}\right)^{-1}}{(2k_1-1)^{c_1}\cdots(2k_d-1)^{c_d}}\nonumber\\
<&\prod_{i=0}^d\frac{\pi|z_i|}{\sin(\pi |z_i|)}t_n^{\star}(\boldsymbol{c})
<\prod_{i=0}^d\frac{\pi u_i}{\sin(\pi u_i)}t^{\star}(\boldsymbol{c}).
\end{align}
Using \eqref{sum-t^*_n,m}, we get
\begin{align*}
		\left|\prod\limits_{k=n+1}^{\infty}\left(1-\frac{z_0^2}{(2k-1)^2}\right)^{-1}-1\right|=\left|\sum_{l=1}^{\infty}t_{\infty,n+1}^{\star}(\{2\}^l)z_0^{2l}\right|<\sum_{l=1}^{\infty}t_{\infty,n+1}^{\star}(\{2\}^l).
\end{align*}
Note that
\begin{align*}
t_{\infty,n+1}^{\star}(\{2\}^l)&=\sum_{k_1\geq\cdots\geq k_l\geq n+1}\frac{1}{(2k_1-1)^2\cdots (2k_l-1)^2}\\
&<\left(\sum_{k=n+1}^{\infty}\frac{1}{(2k-1)^2}\right)^l<\left(\int_{n}^{\infty}\frac{\mathrm{d}x}{(2x-1)^2}\right)^l=\left(\frac{1}{2(2n-1)}\right)^l,
\end{align*}
then we have
\begin{align}\label{1/4n-3}
\left|\prod\limits_{k=n+1}^{\infty}\left(1-\frac{z_0^2}{(2k-1)^2}\right)^{-1}-1\right|<\sum_{l=1}^{\infty}\left(\frac{1}{2(2n-1)}\right)^l=\frac{1}{4n-3}.
\end{align}
Using \eqref{G_n^inf-G_n}, \eqref{G_n-upper-bound} and \eqref{1/4n-3}, we get
\begin{align*}
\left|G_{n}^{\infty}(\boldsymbol{c}; \boldsymbol{z})-G_{n}(\boldsymbol{c}; \boldsymbol{z})\right|
<\frac{1}{4n-3}\prod_{i=0}^d\frac{\pi u_i}{\sin(\pi u_i)}t^{\star}(\boldsymbol{c})\rightarrow0,
\end{align*}
as $n\rightarrow\infty$ on the closed region $E$. Similarly, we have
\begin{align*}
&\left|G(\boldsymbol{c}; \boldsymbol{z})-G_{n}^{\infty}(\boldsymbol{c}; \boldsymbol{z})\right|\\
\leq&\sum_{k_1\geq\cdots\geq k_d\geq1\atop k_1>n}\frac{\prod\limits_{k=k_1}^{\infty}\left(1-\frac{|z_0|^2}{(2k-1)^2}\right)^{-1}\prod\limits_{k=k_2}^{k_1}\left(1-\frac{|z_1|^2}{(2k-1)^2}\right)^{-1}\prod\limits_{k=1}^{k_d}\left(1-\frac{|z_d|^2}{(2k-1)^2}\right)^{-1}}{(2k_1-1)^{c_1}\cdots(2k_d-1)^{c_d}}\nonumber\\
<&\prod_{i=0}^d\frac{\pi u_i}{\sin(\pi u_i)}\left(t^{\star}(\boldsymbol{c})-t_n^{\star}(\boldsymbol{c})\right)\rightarrow0,
\end{align*}
as $n\rightarrow\infty$ on the closed region $E$. Therefore, we finish the proof.
\qed
	
Inspired by \cite[Lemma 4.2]{Linebarger-Zhao}, we get the following lemma.

\begin{lem}\label{limit-W_k-lemma}
Let $a, b, C \in\mathbb{R}$ with $b> 1, C>0$, and let $W_k$ be real numbers satisfying
$$|W_k|< \frac{C\log^a(2k+1)}{k^b}$$
for any $k\in\mathbb{N}$. Then we have
\begin{align*}
\lim_{n\rightarrow\infty}\sum_{k=1}^{n}|W_k|\left(\frac{2}{\pi}-\frac{n}{2^{4n-2}}\binom{2n}{n}\binom{2n-1}{n-k}\right)=0.
\end{align*}
\end{lem}

\proof
In\cite{Wallis1656}, J. Wallis presented the following formula
\begin{align*}
\lim_{n\rightarrow\infty}\left[\frac{(2n)!!}{(2n-1)!!}\right]^2\frac{1}{2n+1}=\frac{\pi}{2},
\end{align*}
and a proof of the Wallis formula relies on evaluation of integrals of powers of $\sin x$, that is, let $I_n=\int_{0}^{\frac{\pi}{2}}\sin^nx\mathrm{d}x$, by repeated partial integration and the value range of $\sin x$, we get
\begin{align*}
I_{n}=\frac{n-1}{n}I_{n-2} \quad \text{and}\quad
I_{2n+1}<I_{2n}<I_{2n-1}
\end{align*}
respectively. Hence we obtain
\begin{align}\label{inequation-double factorial}
\frac{(2n)!!}{(2n+1)!!}<\frac{(2n-1)!!}{(2n)!!}\frac{\pi}{2}<\frac{(2n-2)!!}{(2n-1)!!}.
\end{align}
Using the squeeze theorem, we can deduce the Wallis formula.

As
\begin{align*}
\frac{2}{\pi}-\frac{n}{2^{4n-2}}\binom{2n}{n}\binom{2n-1}{n-k}&=\frac{2}{\pi}-\frac{[(2n-1)!!]^2}{(2n-2k)!!(2n+2k-2)!!}\\
&\geq\frac{2}{\pi}-\frac{[(2n-1)!!]^2}{(2n-2)!!(2n)!!},
\end{align*}
we use \eqref{inequation-double factorial} to get
$$\frac{2}{\pi}-\frac{n}{2^{4n-2}}\binom{2n}{n}\binom{2n-1}{n-k}>0.$$
On the other hand, by \eqref{inequation-double factorial}, we have
\begin{align}\label{upperbound-2/pi-combination}
\frac{2}{\pi}-\frac{n}{2^{4n-2}}\binom{2n}{n}\binom{2n-1}{n-k}
&<\frac{2}{\pi}-\frac{2}{\pi}\frac{[(2n)!!]^2}{(2n+1)(2n-2k)!!(2n+2k-2)!!}.
\end{align}
The right-hand side of \eqref{upperbound-2/pi-combination} is $\frac{2}{\pi}-\frac{2}{\pi}\frac{2(n+k)}{2n+1}\frac{\binom{n}{k}}{\binom{n+k}{k}}$, which implies that
\begin{align*}
\frac{2}{\pi}-\frac{n}{2^{4n-2}}\binom{2n}{n}\binom{2n-1}{n-k}<\frac{2}{\pi}\left(1-\frac{\binom{n}{k}}{\binom{n+k}{k}}\right).
\end{align*}
By the proof of \cite[Lemma 4.2]{Linebarger-Zhao}, we have
\begin{align*}
1-\frac{\binom{n}{k}}{\binom{n+k}{k}}\leq\frac{2k^2}{n}.
\end{align*}
Hence we get
\begin{align*}
\frac{2}{\pi}-\frac{n}{2^{4n-2}}\binom{2n}{n}\binom{2n-1}{n-k}<\frac{4k^2}{\pi n}.
\end{align*}
The subsequent proof is basically the same as the proof of \cite[Lemma 4.2]{Linebarger-Zhao}. Let $\lambda=\min\{\frac{b-1}{2}, \frac{1}{2}\}$, we deduce that
\begin{align*}
&\sum_{k=1}^{n}|W_k|\left(\frac{2}{\pi}-\frac{n}{2^{4n-2}}\binom{2n}{n}\binom{2n-1}{n-k}\right)\\
\leq&\frac{4C}{\pi n}\sum_{k=1}^{\lfloor n^{\lambda}\rfloor}\frac{\log^a(2k+1)}{k^{b-2}}+\frac{2C}{\pi}\sum_{k=\lfloor n^{\lambda}\rfloor+1}^{\infty}\frac{\log^a(2k+1)}{k^{b}}.
\end{align*}
Since $b>1$, the last term becomes $0$ when $n$ tends to infinity. If $b\geq2$, then
\begin{align*}
\frac{4C}{\pi n}\sum_{k=1}^{\lfloor n^{\lambda}\rfloor}\frac{\log^a(2k+1)}{k^{b-2}}\leq\frac{4C}{\pi n}\sum_{k=1}^{\sqrt{n}}\log^a(2k+1)\leq\frac{4C\log^a(2\sqrt{n}+1)}{\pi\sqrt{n}}\rightarrow0
\end{align*}
as $n\rightarrow\infty$. If $1<b<2$, then $\lambda=\frac{b-1}{2}$. Since $\lambda(3-b)=\frac{1}{2}(b-1)(3-b)\leq\frac{1}{2}$, we have
\begin{align*}
\frac{4C}{\pi n}\sum_{k=1}^{\lfloor n^{\lambda}\rfloor}\frac{\log^a(2k+1)}{k^{b-2}}
\leq\frac{4Cn^{\lambda(3-b)}\log^a(2n^{\lambda}+1)}{\pi n}
\leq\frac{4C\log^a(2\sqrt{n}+1)}{\pi\sqrt{n}}
\rightarrow0\quad
\end{align*}
as $n\rightarrow\infty$.  Therefore, we prove the lemma.
\qed
	
\begin{lem}\label{widetlide-G_k_0-bound-lemma}
For $k_0\in\mathbb{N}$, $d\in\mathbb{N}_0$,  $\boldsymbol{c}=(c_1,\ldots,c_d)\in(\mathbb{N}\setminus\{2\})^d$ with $c_1\geq3$, $\boldsymbol{z}=(z_0,z_1,\ldots,z_d)\in\mathbb{C}^{d+1}$
with $|z_j|<1, j=0,1,\ldots,d$, let
\begin{align}\label{Definition of widetlide-G_k_0}
\widetilde{G}_{k_0}(\boldsymbol{c};\boldsymbol{z})=\sum_{k_0\geq k_1\geq\cdots\geq k_d\geq1}\prod_{i=0}^{d}\frac{(-1)^{k_i\delta_i}(2k_i-1)^{\delta_i-1}}{(2k_i-1)^2-z_i^2}V_{k_{i-1},k_i}^{\#}(\{1\}^{c_i-3}),
\end{align}
where $k_{-1}=0$, $\delta_i=\delta(c_i)+\delta(c_{i+1})$ with $c_0=1$ and $c_{d+1}=0$. Then there exist $a, b, C \in\mathbb{R}$ with $b> 1$ and $C>0$ such that for all $k_0\in\mathbb{N}$,
$$\left|\widetilde{G}_{k_0}(\boldsymbol{c};\boldsymbol{z})\right|< \frac{C\log^a(2k_0+1)}{k_0^b}.$$
\end{lem}

\proof
For $i=1,2,\ldots,d$, if $c_i\in\{1, 3\}$,  we have
\begin{align*}
V_{k_{i-1},k_i}^{\#}(\{1\}^{c_i-3})=\left(2\cdot\frac{2k_{i-1}-1}{2k_i-1}\right)^{\triangle(k_{i-1},k_i)}\leq2\cdot\frac{2k_{i-1}-1}{2k_i-1},
\end{align*}
and if $c_i>3$, we have
\begin{align*}
V_{k_{i-1},k_i}^{\#}(\{1\}^{c_i-3})&=\frac{2k_{i-1}-1}{2k_i-1}\sum_{k_{i-1}\geq l_1\geq\cdots\geq l_{c_i-3}\geq k_i}\frac{2^{\triangle(k_{i-1},l_1)+\triangle(l_1,l_2)+\cdots+\triangle(l_{c_i-3},k_i)}}{(2l_1-1)\cdots(2l_{c_i-1}-1)}\\
&<\frac{2k_{i-1}-1}{2k_i-1}\cdot C_{1}\log^{c_i-3}(2k_{i-1}+1),
\end{align*}
where $C_{1}>2$ is some positive constant. Then we can express the above in a unified form
\begin{align*}
V_{k_{i-1},k_i}^{\#}(\{1\}^{c_i-3})
		<\frac{2k_{i-1}-1}{2k_i-1}\cdot C_{1}\log^{c_i-1}(2k_{0}+1).
\end{align*}
Notice that $\delta_0=1$, $\delta_d\in\{2, 3\}$ and $\delta_i\in\{0,1,2\}$ for $i=1,2,\ldots,d-1$. So we find
\begin{align*}
\left|\widetilde{G}_{k_0}(\boldsymbol{c};\boldsymbol{z})\right|&=\frac{2}{(2k_0-1)|(2k_0-1)^2-z_0^2|}\\
&\qquad\times\left|\sum_{k_0\geq k_1\geq\cdots\geq k_d\geq1}\prod_{i=1}^{d}\frac{(-1)^{k_i\delta_i}(2k_i-1)^{\delta_i-1}}{(2k_i-1)^2-z_i^2}V_{k_{i-1},k_i}^{\#}(\{1\}^{c_i-3})\right|\\
&<C_{2}\cdot\frac{\log^{c_1+\cdots+c_d-d}(2k_0+1)}{(2k_0-1)^2-|z_0|^2}\sum_{k_0\geq k_1\geq\cdots\geq k_d\geq1}\prod_{i=1}^{d}\frac{2k_i-1}{(2k_i-1)^2-|z_i|^2},
\end{align*}
where $C_2$ is a positive constant. Let $|z_{\max}|=\max\{|z_1|,\ldots,|z_d|\}$, then
\begin{align*}
\left|\widetilde{G}_{k_0}(\boldsymbol{c};\boldsymbol{z})\right|
&<C_{3}\cdot\frac{\log^{c_1+\cdots+c_d-d}(2k_0+1)}{(2k_0-1)^2-|z_0|^2}\left(\sum_{k=1}^{k_0}\frac{1}{2k-1-|z_{\max}|^2}\right)^d\\
&<C_{4}\cdot\frac{\log^{c_1+\cdots+c_d}(2k_0+1)}{(2k_0-1)^2-|z_{0}|^2}\\
&<C_{5}\cdot\frac{\log^{c_1+\cdots+c_d}(2k_0+1)}{k_0^2},
\end{align*}
where $C_{3}, C_{4}$ and $C_{5}$ are positive constants. Therefore by setting $a=c_1+\cdots+c_d$, $b=2$ and $C=C_5$, we conclude the result.
\qed
	
\subsection{Proofs of the theorems}\label{SubSec:ProfThms3-1,3,4}

We now prove Theorems \ref{G-main-thm}, \ref{t^star-theorem} and \ref{a_u>=1-theorem}.

\noindent \textbf{Proof of Theorem \ref{G-main-thm}}
By Lemma \ref{limit-G_n=G-lemma} and Theorem \ref{G_n-main-thm}, we deduce that
$$G(\boldsymbol{c};\boldsymbol{z})=\lim_{n\rightarrow\infty}G_n(\boldsymbol{c};\boldsymbol{z})
=\lim_{n\rightarrow\infty}\sum_{k_0=1}^{n}\frac{n}{2^{4n-2}}\binom{2n}{n}\binom{2n-1}{n-k_0}\widetilde{G}_{k_0}(\boldsymbol{c};\boldsymbol{z}),$$
where $\widetilde{G}_{k_0}(\boldsymbol{c};\boldsymbol{z})$ is defined in \eqref{Definition of widetlide-G_k_0}. Then using Lemmas \ref{limit-W_k-lemma} and \ref{widetlide-G_k_0-bound-lemma}, we find
$$G(\boldsymbol{c};\boldsymbol{z})=\lim_{n\rightarrow\infty}\sum_{k_0=1}^{n}\frac{2}{\pi}\widetilde{G}_{k_0}(\boldsymbol{c};\boldsymbol{z})
=\frac{2}{\pi}\sum_{k_0=1}^{\infty}\widetilde{G}_{k_0}(\boldsymbol{c};\boldsymbol{z}),$$
which concludes Theorem \ref{G-main-thm}.
\qed
	
\noindent \textbf{Proof of Theorem \ref{t^star-theorem}}
Set $\boldsymbol{a}=(a_0,a_1,$\ldots$,a_d)$ and
\begin{align*}
\widetilde{D}_{k_0}(\boldsymbol{c}; \boldsymbol{a})=\sum_{k_0\geq k_1\geq\cdots\geq k_d\geq1}\prod_{i=0}^{d}\frac{(-1)^{k_i\delta_i}}{(2k_i-1)^{2a_i-\delta_i+3}}V_{k_{i-1},k_i}^{\#}(\{1\}^{c_i-3}).
\end{align*}
By the proof of Lemma \ref{widetlide-G_k_0-bound-lemma}, it is easy to deduce that there exist $a, b, C \in\mathbb{R}$ with $b> 1$ and $C>0$ such that for all $k_0\in\mathbb{N}$,
\begin{align}\label{widetilde-D-k_0}
\left|\widetilde{D}_{k_0}(\boldsymbol{c}; \boldsymbol{a})\right|< \frac{C\log^a(2k_0+1)}{k_0^b}.
\end{align}
Then taking the limit $n\rightarrow\infty$ in Corollary \ref{t_n^star-corollary} and applying \eqref{widetilde-D-k_0} and Lemma \ref{limit-W_k-lemma}, we get the desired result.
\qed
	
\noindent \textbf{Proof of Theorem \ref{a_u>=1-theorem}}
Multiplying \eqref{t^star-formula} by $z_0^{a_0}z_1^{a_1}\cdots z_d^{a_d}$ and summing it over the corresponding set of integers $a_0, a_1,\ldots,a_d$ with $a_0, a_1,\ldots,a_d\geq0$ and $a_u\geq 1$,  the result follows easily.
\qed

%%%%%%%%%%%%%%%%%%%%%%%%%%%%%%%%%%%%%%%%%%%%%%%%%%%%%%%%%%%%%%%%%%%%%%%%%%%%%%%%%%%%%%%%%%%%%%%%%%%%%%%%%%%%%%%%%

\section{Evaluations for multiple $t$-star values}\label{Sec: Evaluations for MtSVs}

As applications of Theorem \ref{t^star-theorem}, we obtain some evaluations of multiple $t$-star values $t^{\star}(\{2\}^{a_0},c_1,\{2\}^{a_1},\ldots,c_d,\{2\}^{a_d})$ with $d=0, 1, 2$ and $c_i\in\{1,3\}$. For more general index, we deduce a connection between multiple $t$-star values and the weighted sum formula for alternating multiple $t$-values. To save space, for an alternating multiple $t$-value $t(s_1,\ldots,s_r;\sigma_1,\ldots,\sigma_r)$, we may put a bar on the top of $s_i$ if $\sigma_i=-1$. For example,
\begin{align*}
&t(\overline{s_1},s_2)=t(s_1, s_2; -1, 1)=\sum_{k_1>k_2\geq1}\frac{(-1)^{k_1}}{(2k_1-1)^{s_1}(2k_2-1)^{s_2}},\\
&t(s_1,\overline{s_2},\overline{s_3})=t(s_1, s_2, s_3; 1, -1, -1)=\sum_{k_1> k_2> k_3\geq1}\frac{(-1)^{k_2}(-1)^{k_3}}{(2k_1-1)^{s_1}(2k_2-1)^{s_2}(2k_3-1)^{s_3}}.
\end{align*}

We remark that the following result has appeared in \cite{Chung2019,Hoffman2019}.

\begin{thm}\label{t^star-2^a-thm}
For any $a\in\mathbb{N}_0$, we have
\begin{align}\label{t^star-2^a-formula}
t^{\star}(\{2\}^{a})=-\frac{4}{\pi}t(\overline{2a+1})=\frac{(-1)^{a}\pi^{2a}E_{2a}}{4^{a}(2a)!},
\end{align}
where $E_n$ are the Euler numbers.
\end{thm}

\proof
Setting $d=0$ in Theorem \ref{t^star-theorem}, we get
\begin{align*}
t^{\star}(\{2\}^{a})=\frac{2}{\pi}\sum_{k=1}^{\infty}\frac{(-1)^k}{(2k-1)^{2a}}\cdot\frac{-2}{2k-1}=-\frac{4}{\pi}t(\overline{2a+1})=\frac{4}{\pi}\beta(2a+1),
\end{align*}
where $\beta(z)$ is the Dirichlet beta function defined by the sum
\begin{align*}
\beta(z)=\sum\limits_{k=0}^\infty\frac{(-1)^{k}}{(2k+1)^z}.
\end{align*}
From \cite{Euler}, we get Euler's identity
\begin{align*}
\beta(2a+1)=\frac{(-1)^{a}\pi^{2a+1}E_{2a}}{4^{a+1}(2a)!},
\end{align*}
which finishes the proof.
\qed
	
\begin{thm}\label{t^star-2-3-2-thm}
For any $a, b\in\mathbb{N}_0$, we have
\begin{align}\label{t^star-2-3-2-formula}
&t^{\star}(\{2\}^{a}, 3, \{2\}^{b})=-\frac{4}{\pi}t(\overline{2a+2b+4})-\frac{8}{\pi}t(\overline{2a+2}, 2b+2)\nonumber\\
=&-\frac{4}{\pi}\sum_{r=1}^{a+b+1}2^{-2r}\left[\binom{2r}{2a+1}(1-2^{-2r})+\binom{2r}{2b+1}\right]\zeta(2r+1)t(\overline{2a+2b+3-2r}).
\end{align}
\end{thm}

\proof
Setting $d=1$, $c_1=3$ in Theorem \ref{t^star-theorem}, we obtain
\begin{align*}
t^{\star}(\{2\}^{a}, 3, \{2\}^{b})
=&-\frac{4}{\pi}\sum_{k_0\geq k_1\geq1}\frac{(-1)^{k_0}}{(2k_0-1)^{2a+3}(2k_1-1)^{2b+1}}\left(2\cdot\frac{2k_0-1}{2k_1-1}\right)^{\triangle(k_0,k_1)}\\
=&-\frac{4}{\pi}\sum_{k=1}^{\infty}\frac{(-1)^{k}}{(2k-1)^{2a+2b+4}}-\frac{8}{\pi}\sum_{k_0>k_1\geq1}\frac{(-1)^{k_0}}{(2k_0-1)^{2a+2}(2k_1-1)^{2b+2}}\\
=&-\frac{4}{\pi}t(\overline{2a+2b+4})-\frac{8}{\pi}t(\overline{2a+2}, 2b+2).
\end{align*}
Using \cite[(2.4)]{Quan2020} or \cite[(4.6)]{Xu-Yan}, we have
\begin{align*}
t(\overline{2a+2}, 2b+2)=&-\frac{1}{2}t(\overline{2a+2b+4})\\
&+\sum_{k=0}^{b}\dbinom{2a+2b+2-2k}{2a+1}\frac{\overline{\zeta}(2a+2b+3-2k)}{2^{2a+2b+3-2k}}t(\overline{2k+1})\\
&+\sum_{l=0}^{a}\dbinom{2a+2b+2-2l}{2b+1}\frac{\zeta(2a+2b+3-2l)}{2^{2a+2b+3-2l}}t(\overline{2l+1}),
\end{align*}
where $\overline{\zeta}(z)=\sum\limits_{n=1}^{\infty}\frac{(-1)^{n-1}}{n^z}$. Then
\begin{align*}
t^{\star}(\{2\}^{a}, 3, \{2\}^{b})
=&-\frac{8}{\pi}\sum_{k=0}^{b}\dbinom{2a+2b+2-2k}{2a+1}\frac{\overline{\zeta}(2a+2b+3-2k)}{2^{2a+2b+3-2k}}t(\overline{2k+1})\\
&-\frac{8}{\pi}\sum_{l=0}^{a}\dbinom{2a+2b+2-2l}{2b+1}\frac{\zeta(2a+2b+3-2l)}{2^{2a+2b+3-2l}}t(\overline{2l+1}).
\end{align*}
Note that $\overline{\zeta}(k)=(1-2^{1-k})\zeta(k)$ for $k>1$. Hence we deduce that
\begin{align*}
&t^{\star}(\{2\}^{a}, 3, \{2\}^{b})\\
=&-\frac{8}{\pi}\sum_{k=0}^{b}\dbinom{2a+2b+2-2k}{2a+1}(1-2^{-(2a+2b+2-2k)})\frac{\zeta(2a+2b+3-2k)}{2^{2a+2b+3-2k}}t(\overline{2k+1})\\
&-\frac{8}{\pi}\sum_{l=0}^{a}\dbinom{2a+2b+2-2l}{2b+1}\frac{\zeta(2a+2b+3-2l)}{2^{2a+2b+3-2l}}t(\overline{2l+1}),
\end{align*}
which easily implies the result.
\qed
	
The evaluation of $t^{\star}(\{2\}^{a}, 3, \{2\}^{b})$ in \cite[Theorem 4.1]{Li-Wang} shows that
\begin{align}\label{t^star-2-3-2-formula-Li-Wang}
&t^{\star}(\{2\}^{a}, 3, \{2\}^{b})\nonumber\\
=&\sum_{r=1}^{a+b+1}2^{-2r}\left[\binom{2r}{2a+1}(1-2^{-2r})+\binom{2r}{2b+1}\right]\zeta(2r+1)t^{\star}(\{2\}^{a+b+1-r}).
\end{align}
Applying \eqref{t^star-2^a-formula}, we get
\begin{align*}
-\frac{4}{\pi}t(\overline{2a+2b+3-2r})=t^{\star}(\{2\}^{a+b+1-r}).
\end{align*}
Therefore, we  imply the equivalence of the evaluation formulas \eqref{t^star-2-3-2-formula} and \eqref{t^star-2-3-2-formula-Li-Wang}.
	
\begin{thm}\label{t^star-2-1-2-thm}
For any $a\in\mathbb{N}$, $b\in\mathbb{N}_0$, we have
\begin{align}\label{t^star-2-1-2-formula}
&t^{\star}(\{2\}^{a}, 1, \{2\}^{b})=-\frac{4}{\pi}t(\overline{2a+2b+2})-\frac{8}{\pi}t(2a+1, \overline{2b+1})\nonumber\\
=&-\frac{4}{\pi}\sum_{r=1}^{a+b}2^{-2r}\left[\binom{2r}{2a}+\binom{2r}{2b}(1-2^{-2r})\right]\zeta(2r+1)t(\overline{2a+2b+1-2r})\nonumber\\
&-\delta_{b,0}\frac{4\log2}{\pi}t(\overline{2a+1}),
\end{align}
where $\delta$ denotes the Kronecker symbol.
\end{thm}

\proof
Let $d=1$, $c_1=1$ in Theorem \ref{t^star-theorem}, we obtain
\begin{align*}
t^{\star}(\{2\}^{a}, 1, \{2\}^{b})
=&-\frac{4}{\pi}\sum_{k_0\geq k_1\geq1}\frac{(-1)^{k_1}}{(2k_0-1)^{2a+2}(2k_1-1)^{2b}}\left(2\cdot\frac{2k_0-1}{2k_1-1}\right)^{\triangle(k_0,k_1)}\\
=&-\frac{4}{\pi}\sum_{k=1}^{\infty}\frac{(-1)^{k}}{(2k-1)^{2a+2b+2}}-\frac{8}{\pi}\sum_{k_0>k_1\geq1}\frac{(-1)^{k_1}}{(2k_0-1)^{2a+1}(2k_1-1)^{2b+1}}\\
=&-\frac{4}{\pi}t(\overline{2a+2b+2})-\frac{8}{\pi}t(2a+1, \overline{2b+1}).
\end{align*}
By \cite[(2.2)]{Quan2020} and notice $\overline{\zeta}(k)=(1-2^{1-k})\zeta(k)$ for $k>1$, we easily get the desired result.
\qed
	
Similarly, using \eqref{t^star-2^a-formula}, we can deduce that the evaluation formula \eqref{t^star-2-1-2-formula} is equivalent to \cite[(4.2)]{Li-Wang}.
	
The following theorem shows that a special weighted sum formula for alternating multiple $t$-values can be expressed by multiple $t$-star values with two blocks of twos.

\begin{thm}\label{t^star-2-c+3-2-thm}
For any $a, b\in\mathbb{N}_0$ and $c\in\mathbb{N}$, we have
\begin{align*}
t^{\star}(\{2\}^{a}, c+3, \{2\}^{b})
=-\frac{2}{\pi}\sum_{r=1}^{c+2}2^{r}\sum_{s_1+\cdots+s_r=c+2\atop s_1,\ldots,s_r\geq1}t(\overline{2a+1+s_1}, s_2,\ldots, s_{r-1}, 2b+1+s_r).
\end{align*}
Here $t(\overline{2a+1+s_1}, s_2,\ldots, s_{r-1}, 2b+1+s_r)=t(\overline{2a+2b+4+c})$ when $r=1$.
\end{thm}

\proof
Setting $d=1$, $c_1=c+3 $ in Theorem \ref{t^star-theorem}, we obtain that
\begin{align*}
&t^{\star}(\{2\}^{a}, c+3, \{2\}^{b})\\
=&-\frac{4}{\pi}\sum_{k_0\geq k_1\geq1}\frac{(-1)^{k_0}}{(2k_0-1)^{2a+3}(2k_1-1)^{2b+1}}V_{k_0,k_1}^{\#}(\{1\}^c)\\
=&-\frac{4}{\pi}\sum_{k_0\geq l_1\geq\cdots\geq l_c\geq k_1\geq1}\frac{(-1)^{k_0}2^{\triangle(k_0, l_1)+\triangle(l_1, l_2)+\cdots+\triangle(l_{c-1},l_c)+\triangle(l_c, k_1)}}{(2k_0-1)^{2a+2}(2l_1-1)\cdots(2l_c-1)(2k_1-1)^{2b+2}}.
\end{align*}
By counting the numbers of $k_0, l_1, \ldots, l_c, k_1$ which are not equal, we have
\begin{align*}
&t^{\star}(\{2\}^{a}, c+3, \{2\}^{b})\\
=&-\frac{4}{\pi}\sum_{k=1}^{\infty}\frac{(-1)^k}{(2k-1)^{2a+2b+4+c}}-\frac{4}{\pi}\sum_{r=2}^{c+2}2^{r-1}\sum_{s_1+\cdots+s_r=c+2\atop s_1,\ldots,s_r\geq1}\\
&\quad\quad\times\sum_{k_1>\cdots>k_r\geq1}\frac{(-1)^{k_1}}{(2k_1-1)^{2a+1+s_1}(2k_2-1)^{s_2}\cdots(2k_{r-1}-1)^{s_{r-1}}(2k_r-1)^{2b+1+s_r}}\\
=&-\frac{2}{\pi}\sum_{r=1}^{c+2}2^{r}\sum_{s_1+\cdots+s_r=c+2\atop s_1,\ldots,s_r\geq1}t(\overline{2a+1+s_1}, s_2,\ldots, s_{r-1}, 2b+1+s_r),
\end{align*}
where if $r=1$,  $t(\overline{2a+1+s_1}, s_2,\ldots, s_{r-1}, 2b+1+s_r)$ is treated as $t(\overline{2a+2b+4+c})$.
\qed
	
For $d=2$, setting $(c_1,c_2)=(3, 1), (3, 3)$ in Theorem \ref{t^star-theorem}, we deduce the following theorem.

\begin{thm}
For any $a, b, c\in\mathbb{N}_0$, we have
\begin{align*}
t^{\star}(\{2\}^{a}, 3, \{2\}^{b}, 1, \{2\}^{c})
=&-\frac{4}{\pi}t(\overline{2a+2b+2c+5})-\frac{8}{\pi}t(\overline{2a+2}, 2b+2c+3)\\
&-\frac{8}{\pi}t(2a+2b+4, \overline{2c+1})-\frac{16}{\pi}t(\overline{2a+2}, \overline{2b+2}, \overline{2c+1})
\end{align*}
and
\begin{align*}
t^{\star}(\{2\}^{a}, 3, \{2\}^{b}, 3, \{2\}^{c})
=&-\frac{4}{\pi}t(\overline{2a+2b+2c+7})-\frac{8}{\pi}t(\overline{2a+2}, 2b+2c+5)\\
&-\frac{8}{\pi}t(\overline{2a+2b+5}, 2c+2)-\frac{16}{\pi}t(\overline{2a+2}, 2b+3, 2c+2).
\end{align*}
\end{thm}
For $(c_1,c_2)=(1, 3), (1,1)$ , we can get similar results as above.
\begin{thm}
For any $a\in\mathbb{N}$, $b, c\in\mathbb{N}_0$, we have
\begin{align*}
t^{\star}(\{2\}^{a}, 1, \{2\}^{b}, 3, \{2\}^{c})
=&-\frac{4}{\pi}t(\overline{2a+2b+2c+5})-\frac{8}{\pi}t(2a+1, \overline{2b+2c+4})\\
&-\frac{8}{\pi}t(\overline{2a+2b+3}, 2c+2)-\frac{16}{\pi}t(2a+1, \overline{2b+2}, 2c+2)
\end{align*}
and
\begin{align*}
t^{\star}(\{2\}^{a}, 1, \{2\}^{b}, 1, \{2\}^{c})
=&-\frac{4}{\pi}t(\overline{2a+2b+2c+3})-\frac{8}{\pi}t(2a+1, \overline{2b+2c+2})\\
&-\frac{8}{\pi}t(2a+2b+2, \overline{2c+1})-\frac{16}{\pi}t(2a+1, 2b+1, \overline{2c+1}).
\end{align*}
\end{thm}

Setting $c_1=\cdots=c_{d}=1$, $a_0=a_d=1$ and $a_1=\cdots=a_{d-1}=0$ in Theorem \ref{t^star-theorem}, we have the following theorem.

\begin{thm}
For any $d\in\mathbb{N}_0$, we have
	\begin{align*}
		t^{\star}(2,\{1\}^d,2)
		=-\frac{2}{\pi}\sum_{r=1}^{d+1}2^r\sum_{s_1+\cdots+s_r=d+1\atop s_1,\ldots,s_r\geq1}t(s_1+2,s_2,\ldots,s_{r-1},\overline{s_r+2}).
	\end{align*}
	Here for $r=1$, the inner sum becomes the single alternating $t$-value $t(\overline{d+5})$.
\end{thm}

Finally, setting $(c_1,\ldots,c_{2d-2},c_{2d-1})=(\{3,1\}^{d-1}, 3)$ and $a_0=a_1=\cdots=a_{2d-1}=a$ with $a\geq0$ in Theorem \ref{t^star-theorem}, we obtain the following result.

\begin{thm}
For any $d\in\mathbb{N}$ and $a\in\mathbb{N}_0$, we have
\begin{align*}
&t^{\star}(\{\{2\}^a,3,\{2\}^a,1\}^{d-1},\{2\}^a,3,\{2\}^a)\\
=&-\frac{2}{\pi}\sum_{r=1}^{2d}2^r\sum_{s_1+\cdots+s_r=2d\atop s_1,\ldots,s_r\geq1}t((2a+2)s_1,\ldots,(2a+2)s_r;  (-1)^{s_1},\ldots,(-1)^{s_{r-1}},(-1)^{s_r-1}).
\end{align*}
Here for $r=1$, the inner sum becomes the single alternating $t$-value $t(\overline{2d(2a+2)})$.
\end{thm}

%%%%%%%%%%%%%%%%%%%%%%%%%%%%%%%%%%%%%%%%%%%%%%%%%%%%%%%%%%%%%%%%%%%%%%%%%%%%%%%%%%%%%%%%%%%%%%%%%%%%%%%%%%%%%%%%

\end{document}